\documentclass[a4paper,11pt]{amsart}

\usepackage[french]{babel}
\usepackage{amsmath, enumerate, amsfonts, amssymb, amsthm}

\newtheorem{defin}{\bf D\'ef\mbox{}inition}[section]
\newtheorem{theo}[defin]{\bf Th\'eor\`eme}
\newtheorem{prop}[defin]{\bf Proposition}
\newtheorem{lem}[defin]{\bf Lemme}
\newtheorem{cor}[defin]{\bf Corollaire}
\newtheorem{propri}[defin]{\bf Propri\'et\'e}
\newtheorem{aff}[defin]{\bf Af\mbox{}f\mbox{}irmation}
\newtheorem{rem}[defin]{\bf Remarque}
\newtheorem{nota}[defin]{\bf Notation}
\newtheorem*{theoBMM}{\bf Th\'eor\`eme BMM}

\newtheorem{theoi}{\bf Th\'eor\`eme}
\newtheorem{propi}{\bf Proposition}
\newtheorem*{ex-di}{\bf Exemples et Discussions}

\newcommand{\dps}{\displaystyle}

\newcommand{\gauche}{\begin{flushleft}\end{flushleft}}

\newcommand{\R}{\mathbb{R}}
\newcommand{\C}{\mathbb{C}}
\newcommand{\N}{\mathbb{N}}

\newcommand{\Q}{\mathbb{Q}}

\newcommand{\D}{\mathcal{D}}
\newcommand{\Dn}{\mathcal{D}_n}
\newcommand{\dx}[1]{\partial _{x_{#1}}}
\newcommand{\ddx}{\partial _x}

\newcommand{\dxsur}[2]{\frac{\partial {#1}}{\partial x_{#2}}}

\newcommand{\B}{\mathcal{B}}
\newcommand{\rhoun}{\rho_{\alpha_1}}
\newcommand{\rhode}{\rho_{\alpha_2}}
\newcommand{\degun}{\mathrm{deg}_{\alpha_1}}
\newcommand{\degde}{\mathrm{deg}_{\alpha_2}}
\newcommand{\fin}{\mathrm{f{}in}}
\newcommand{\ini}{\mathrm{in}}
\newcommand{\inu}{\mathrm{in_{\alpha_1}}}
\newcommand{\ind}{\mathrm{in_{\alpha_2}}}
\newcommand{\cp}{\mathrm{cp}} 
\newcommand{\tp}{\mathrm{tp}}
\newcommand{\DN}{\mathcal{N}} 

\title[Une $b$-fonction remarquable]
{Construction d'un \'el\'ement remarquable de l'id\'eal de Bernstein-Sato
associ\'e \`a deux courbes planes analytiques}

\author{Rouchdi BAHLOUL}
\address{Department of Mathematics, Faculty of Science, Kobe University,
1-1, Rokkodai, Nada-ku, Kobe 657-8501, Japan}
\email{rouchdi@math.kobe-u.ac.jp}

\begin{document}

\begin{abstract}
Etant donn\'ees des d\'eformations $f_1,f_2$ de deux polyn\^omes de
deux variables quasi-homog\`enes pour deux syst\`emes de poids
distincts $\alpha_1, \alpha_2$ satisfaisant \`a des conditions similaires
\`a celles de singularit\'es semi-quasi-homog\`enes pour un poids,
nous donnons, par des m\'ethodes inspir\'ees de celles de
H. Maynadier, une formule explicite d'un polyn\^ome de Bernstein-Sato
faisant intervenir deux formes affines $\rho_{\alpha_i}(f_1) s_1 +
\rho_{\alpha_i}(f_2) s_2 +k$, $i=1,2$. Dans le cas particulier
$(f_1,f_2)=(x_1^a+x_2^b, x_1^c+x_2^d)$, $bc-ad >0$, nous calculons
l'espace $\mathcal{H}_f$ \'etudi\'e r\'ecemment par J. Brian\c{c}on,
Ph. Maisonobe et M. Merle et montrons qu'il est \'egal au lieu des z\'eros
de $s_1 s_2 (ab s_1+ ads_2) (ads_1+ cds_2)$.
\end{abstract}


\maketitle

\section*{Introduction et \'enonc\'e des r\'esultats principaux}

L'\'etude du polyn\^ome de Bernstein (ou $b$-fonction selon M. Sato)
et de son analogue \`a plusieurs fonctions, l'id\'eal de Bernstein-Sato
remonte au d\'ebut des ann\'ees 1970 avec entre autres les contributions
fondamentales de I.~N. Bernstein \cite{bernstein}, J.E. Bj\"ork \cite{bjork},
B.~Malgrange \cite{malgrange} et M.~Kashiwara \cite{kashiwara} pour le
polyn\^ome de Bernstein et de C.~Sabbah \cite{sabbah1, sabbah2} pour
l'id\'eal (ou les id\'eaux) de Bernstein-Sato.

Ces \'etudes ``th\'eoriques'' n'ont en rien \'epuis\'e le sujet. En
effet, aujourd'hui encore, il reste beaucoup de questions ouvertes,
notamment sur des questions de rationalit\'e ou de principalit\'e, mais
aussi concernant le calcul explicite de ces polyn\^omes et de ces id\'eaux.
\`A ce propos, citons les contributions de T.~Yano et M.~Kato qui ont
trait\'e nombre d'exemples. Citons aussi T.~Oaku \cite{oaku} qui, le premier,
donna un algorithme de calcul du polyn\^ome de Bernstein global et local
d'un polyn\^ome, et ce sans aucune hypoth\`ese sur le polyn\^ome
en question; cet algorithme est bas\'e sur la th\'eorie des bases de
Gr\"obner dans les anneaux d'op\'erateurs diff\'erentiels. Signalons
\cite{walther} dans lequel U.~Walther calcule le polyn\^ome de Bernstein
(global) d'un arrangement d'hyperplans g\'en\'erique. Enfin, terminons
cette liste (non exhaustive) par une contribution qui nous concerne
directement ici, \`a savoir l'article de J.~Brian\c{c}on et \emph{al}
\cite{bgmm} dans lequel les auteurs donnent une m\'ethode explicite pour
d\'eterminer le polyn\^ome de Bernstein d'une singularit\'e
semi-quasi-homog\`ene ou non d\'eg\'en\'er\'ee au sens de Kouchnirenko
(voir aussi T.~Torrelli \cite{torrelli}).

Ici nous nous int\'eressons aux cas de plusieurs fonctions.
\`A $p$ germes $f_1,\ldots,f_p \in \C\{x\}=\C\{x_1,\ldots,x_n\}$, on
associe l'id\'eal de Bernstein-Sato $\B$ form\'e des polyn\^omes
$b(s) \in \C[s]=\C[s_1,\ldots,s_p]$ satisfaisant \`a
\[b(s) f_1^{s_1}\cdots f_p^{s_p} \in \Dn[s] F f_1^{s_1}\cdots f_p^{s_p}\]
o\`u $\Dn$ est l'anneau des op\'erateurs diff\'erentiels \`a coefficients
dans $\C\{x\}$ et $F$ est le produit des $f_j$.
D'apr\`es C. Sabbah \cite{sabbah1, sabbah2}, cet id\'eal n'est
pas nul et d'apr\`es J. Brian\c{c}on et H. Maynadier \cite{bmay}, il n'est
pas principal en g\'en\'eral.

H.~Maynadier \cite{maynadierTez, maynadierBull} a \'etudi\'e l'id\'eal de
Bernstein-Sato d'une singularit\'e quasi-homog\`ene ainsi que
semi-quasi-homog\`ene (i.e. l'id\'eal des parties initiales d\'efinit
une intersection compl\`ete \`a singularit\'e isol\'ee)
pour \emph{un} syst\`eme de poids donn\'e en s'inspirant de \cite{bgmm}.
La m\'ethode peut se r\'esumer ainsi~: dans un premier temps, on fait
monter le poids des coefficients des op\'erateurs \`a l'aide de
l'op\'erateur d'Euler associ\'e au syst\`eme de poids et dans un deuxi\`eme
temps on r\'e\'ecrit ces coefficients via des divisions par des id\'eaux
de colongueur finie.

Une question naturelle se pose alors~: que se passe-t-il si on essaie de
g\'en\'eraliser ces m\'ethodes dans le cas o\`u l'on s'autorise des
syst\`emes de poids distincts~? La premi\`ere situation \`a \'etudier est
celle de deux fonctions de deux variables. C'est ce qu'on se propose de
faire ici.\\

{\bf Notations.}
Pour un germe de fonction analytique $u=\sum_{ij} u_{ij}x_1^i x_2^j$,
avec $u_{ij} \in \C$, on d\'efinit son diagramme de Newton $\DN(u)
\subset \N^2$ comme l'ensemble des couples $(i,j)$ tels que $u_{ij}
\ne 0$. \'Etant donn\'e un syst\`eme de poids $\alpha=(b,a)\in (\Q_{>0})^2$,
un polyn\^ome $p=p(x_1,x_2)$ est dit $\alpha$-homog\`ene de degr\'e
$\deg_\alpha(p)=d$ si pour tout $(i,j)\in \DN(p)$, $bi+aj=d$. Pour un
polyn\^ome quelconque $p$, le degr\'e par rapport \`a $\alpha$ est
$\deg_{\alpha}(p)=\max\{ b i+a j/\, (i,j)\in \DN(p)\}$. 
On d\'efinit le poids (pour le syst\`eme $\alpha$) d'une s\'erie $u$ comme
$\rho_\alpha(u)= \min\{ b i+a j/\, (i,j)\in \DN(u)\}$. On note $\ini_\alpha(u)$
la partie initiale de $u$, i.e. sa partie $\alpha$-homog\`ene de degr\'e
$\rho_\alpha(u)$. De plus, pour un id\'eal $I$ de $\C\{x_1,x_2\}$, on note
$\ini_\alpha(I)$ l'id\'eal de $\C\{x_1,x_2\}$ engendr\'e par les parties
initiales des \'el\'ements de $I$. Enfin, pour $u_1, \ldots, u_r \in
\C\{x_1,x_2\}$, on note $\langle u_1, \ldots, u_r \rangle$ l'id\'eal qu'ils
engendrent.\\

On se donne deux fonctions analytiques complexes $f_1,f_2$ de deux
variables $x_1,x_2$, deux entiers positifs non nuls $a$ et $d$ et deux
rationnels $b$ et $c$ strictement sup\'erieurs \`a $1$ tels que les
deux syst\`emes de poids $\alpha_1=(b,a)$, $\alpha_2=(d,c)$ soient non
colin\'eaires. On suppose~:
\begin{enumerate}
\item $x_1^a$ (resp. $x_2^d$) est un terme de $f_1$ (resp. $f_2$),\\
i.e. $(a,0) \in \DN(f_1)$ et $(0,d) \in \DN(f_2)$,
\item $\rhoun(f_1)=\degun(x_1^a)=ab$,
\item $\rhode(f_2)=\degde(x_2^d)=cd$,
\item $bc > ad$ (i.e. $\inu(f_1)$ est ``plus pentue'' que $\ind(f_2)$).
\end{enumerate}
On impose les hypoth\`eses suppl\'ementaires suivantes.  Si on note
$J$ le d\'eterminant jacobien de $(f_1,f_2)$, alors on demande~:
\begin{enumerate}
\item[(5)] Les id\'eaux $\inu (\langle f_1, J \rangle)$ et $\langle
  \inu(f_1), \inu(J) \rangle$ sont \'egaux et de colongueur finie dans
  $\C\{x_1,x_2\}$.
\item[(6)] Les id\'eaux $\ind (\langle f_2, J \rangle)$ et $\langle
  \ind(f_2), \ind(J) \rangle$ sont \'egaux et de colongueur finie dans
  $\C\{x_1,x_2\}$.\\
\end{enumerate}
Notons $\D=\D_2$ l'anneau des op\'erateurs diff\'erentiels \`a
coefficients dans $\C\{x\}=\C\{x_1,x_2\}$. Le r\'esultat principal de
cette note est le suivant.

\begin{theoi}\label{theoi}
\gauche
\begin{itemize}
\item
Si $a=1$ et $d=1$ alors l'id\'eal de Bernstein-Sato $\B$ est engendr\'e
par $(s_1+1)(s_2+1)$.
\item
Si $a\ge 2$ ou $d\ge 2$~:\\
Posons $N_1=2ab+ad-2a-2b$, $N_2=2cd+ad-2c-2d$,
\[W_1=\{\degun(m)/ \, m \text{ mon\^ome et } \degun(m)\le N_1+
\rhoun(f_2)\},\]
\[W_2=\{\degde(m)/ \, m \text{ mon\^ome et } \degde(m)\le N_2+
\rhode(f_1)\}.\]
Le polyn\^ome
\[b(s_1,s_2)=(s_1+1) (s_2+1) \prod_{\rho_1 \in W_1} (abs_1
+ads_2+ a+b+ \rho_1) \cdot \prod_{\rho_2 \in W_2} (ads_1 +cds_2+ c+d+
\rho_2)\]
appartient \`a l'id\'eal $\B$.
\end{itemize}
\end{theoi}

La construction que nous mettons en place s'applique \`a la situation
suivante~:

\begin{propi}\label{propi1}
Soient $a,b,c,d$ des entiers non nuls tels que $bc> ad$.  Soient $g_1,
g_2 \in \C\{x\}$. Si $g_1$ (resp. $g_2$) est non nul, on suppose
$\rhoun(g_1) > ab$ (resp. $\rhode(g_2) > cd$).
Posons $f_1=x_1^a+x_2^b+g_1$ et
$f_2=x_1^c+x_2^d+g_2$. Alors $(f_1,f_2)$ satisfait aux hypoth\`eses
(1),\ldots,(6) et le th\'eor\`eme pr\'ec\'edent s'applique.
\end{propi}

\begin{ex-di}
\rm
\begin{itemize}
\item[(1)]
Consid\'erons $f=(x_1, x_1^3+x_2^2)$. Il est facile de voir que si
l'on pose $\alpha_1=(b,1)$ et $\alpha_2=(2,3)$ alors les hypoth\`eses
du th\'eor\`eme sont satisfaites pour n'importe quel entier $b\ge 2$.
On remarque que $N_2=4$ et $W_2=\{0,2,3,4,5,6\}$ quelque soit $b$,
alors que $N_1$ et $W_1$ d\'ependent de $b$.
Ainsi le polyn\^ome suivant $(s_1+1)(s_2+1) p_b(s_1,s_2) q(s_1,s_2)$,
o\`u $p_b(s_1,s_2)=\prod_{\rho_1 \in W_1} (bs_1 +2s_2+ 1+b+ \rho_1)$
et $q(s_1,s_2)=\prod_{\rho_2 \in W_2} (2s_1 +6s_2+ 5+ \rho_2)$,
est dans $\B$, ce quelque soit $b\ge 2$. Ceci entraine que
$b_1(s):=(s_1+1)(s_2+1)\prod_{\rho_2 =0,2,3,4,5,6} (2s_1 +6s_2+ 5+ \rho_2)$
est lui-m\^eme dans $\B$.
\item[(2)]
H. Maynadier \cite{maynadierBull} calcule l'id\'eal de Bernstein-Sato
associ\'e \`a deux polyn\^omes $f_1,f_2$ de deux variables,
$\alpha$-homog\`enes pour un syst\`eme de poids $\alpha$ donn\'e,
dans le cas o\`u les
$f_j$ sont \`a singularit\'e isol\'ee et le morphisme $f=(f_1,f_2)$
d\'efinit l'origine. Elle montre alors que l'id\'eal en question
est principal et en calcule le g\'en\'erateur. Dans l'exemple
pr\'ec\'edent (avec $\alpha=(2,3)$), elle obtient comme g\'en\'erateur
$b_0(s):=(s_1+1)(s_2+1)\prod_{\rho_2 =0,2,4,6} (2s_1 +6s_2+ 5+ \rho_2)$.
Ce dernier divise $b_1$. Cela montre en passant que les $W_j$
ne sont pas optimaux en g\'en\'eral.
\item[(3)]
Dans \cite{maynadierBull}, H. Maynadier traite aussi les morphismes
$f=(f_1,\ldots,f_p)$ tels que toute famille extraite d\'efinit une
singularit\'e semi-quasi-homog\`ene pour $\alpha$ (i.e. les formes
initiales de ses fonctions forment une intersection compl\`ete \`a
singularit\'e isol\'ee). Dans ce cas, elle montre qu'un polyn\^ome
de la forme $\prod_{j=1}^p (s_j+1) \prod_{\rho\in W} (\sum_{j=1}^p
\rho_\alpha(f_j) s_j +\rho)$, $W$ \'etant un ensemble d'entiers,
appartient \`a $\B$.

Maintenant, consid\'erons $f=(x_1^2+x_2^3, x_1^3+x_2^2)$.
Signalons qu'on ne connait
toujours pas l'id\'eal de Bernstein-Sato de $f$ m\^eme dans le cas
global, cf. \cite{bahloulJ} et \cite{ucha-castro}.
On peut facilement voir que cet exemple ne rentre dans aucune des
situations \'etudi\'ees dans \cite{maynadierBull}.
La proposition pr\'ec\'edente nous donne cependant un polyn\^ome dans
$\B$, \`a savoir
\[(s_1+1)(s_2+1)\prod_{\rho}(6s_1+4s_2+5+\rho)(4s_1+6s_2+5+\rho)\]
avec $\rho$ qui parcourt $W_1=W_2=\{0,2,3,\ldots,10\}$.
Au vu de l'exemple (1), on peut penser qu'un produit sur un sous-ensemble
strict de $W_1$ suffise
. Par ailleurs on peut se demander si
les deux formes affines qui apparaissent sont obligatoires.
La r\'eponse est oui comme on va le voir plus bas et cela montre
que la situation examin\'ee dans ce papier est bien diff\'erente
de celles trait\'es par H. Maynadier.
\end{itemize}
\end{ex-di}

De fa\c{c}on g\'en\'erale : soient $f_1, \ldots, f_p \in \C\{x\}$ des
fonctions analytiques au voisinage de $0\in \C^n$. Au voisinage de $x=0$,
soit $A$ le sous-ensemble de $T^* \C^n \times \C^p$ d\'efini par~:
\[A=\big\{ (x, \lambda_1 df_1(x)+\cdots+ \lambda_p df_p(x),
\lambda_1 f_1(x), \ldots, \lambda_p f_p(x) )\big\},\]
o\`u $(\lambda_1, \ldots, \lambda_p)$ d\'ecrit $\C^p$.  Consid\'erons
l'espace $W_f^\#$ introduit par M.~Kashiwara et T.~Kawai \cite{kk}
d\'efini comme l'adh\'erence de $A$ dans $T^* \C^n \times \C^p$; il s'agit
de la vari\'et\'e caract\'eristique de $\D_n[s_1,\ldots,s_p] \cdot f_1^{s_1}
\cdots f_p^{s_p}$ vu comme $\D_n[s_1,\ldots,s_p]$-module (voir \cite{kk},
Th. 2). Maintenant, suivant J. Brian\c{c}on, Ph. Maisonobe et M. Merle
\cite{bmm1, bmm2}, notons $\pi_2: T^* \C^n \times \C^p \to \C^p$ la
projection canonique et consid\'erons $\mathcal{H}_f \subset \C^p$
d\'efini par~:
\[\mathcal{H}_f= \pi_2(W_f^\# \cap F^{-1}(0)),\]
o\`u l'on \'etend canoniquement $F$ \`a $T^* \C^n \times \C^p$ par
$F(x,\xi,s)=F(x)$.\\
Pour un polyn\^ome $c \in \C[s_1,\ldots, s_p]$, on note $\fin(c)$ sa
partie homog\`ene de plus haut degr\'e (sa partie f{}inale) et pour un
id\'eal $J \subset \C[s_1,\ldots, s_p]$, $\fin(J)$ d\'esigne l'id\'eal
engendr\'e par l'ensemble des parties f{}inales des $c \in J$.
\begin{theoBMM}[\cite{bmm1, bmm2}]
\
\begin{enumerate}
\item $\mathcal{H}_f$ est une r\'eunion d'hyperplans vectoriels.
\item $\mathcal{H}_f$ est contenue dans le lieu des z\'eros de
$\fin(\B)$.
\item Si $p=2$~: il existe $b \in \B$ tel que le lieu des z\'eros de
$\fin(b)$ \'egale $\mathcal{H}_f$. Par cons\'equent $\mathcal{H}_f$
est \'egal au lieu des z\'eros de $\fin(\B)$.
\end{enumerate}
\end{theoBMM}

Revenons \`a notre situation. D'apr\`es ce th\'eor\`eme et le
th\'eor\`eme \ref{theoi},
\[\mathcal{H}_f \subset V(s_1 s_2 (ab s_1 + ad s_2) (ad s_1+ cd s_2)).\]

\begin{propi}\label{propi2}
Dans le cas o\`u $(f_1,f_2)= (x_1^a+x_2^b, x_1^c+x_2^d)$ et $bc > ad$,
$\mathcal{H}_f$ est \'egal \`a $V(s_1 s_2 (ab s_1 + ad s_2) (ad s_1+
cd s_2))$.
\end{propi}
Ainsi dans ce cas, le polyn\^ome $b$ construit dans le
th\'eor\`eme \ref{theoi} v\'erifie $V(\fin(b))=\mathcal{H}_f$, on
retrouve donc le r\'esultat (3) du th\'eor\`eme pr\'ec\'edent dans ce
cas bien particulier. Cela justifie l'adjectif ``remarquable''
employ\'e dans le titre.

Comme cons\'equence du th\'eor\`eme pr\'ec\'edent, tout polyn\^ome de
Bernstein-Sato $b$ associ\'e \`a $(f_1,f_2)= (x_1^a+x_2^b,
x_1^c+x_2^d)$ ($bc > ad$) a comme facteurs obligatoires de sa partie
f{}inale les formes lin\'eaires $(ab s_1 + ad s_2)$ et $(ad s_1+ cd s_2)$.
On peut conjecturer que l'\'egalit\'e de la proposition \ref{propi2}
soit encore vraie sous les hypoth\`eses plus g\'en\'erales du
th\'eor\`eme \ref{theoi} ou au moins celles de la proposition
\ref{propi1}.

Signalons qu'en g\'en\'eral, nous ne savons pas determiner
les formes lin\'eaires obligatoires. Nous savons qu'elles
sont contenues dans le $1$-squelette du $V$-\'eventail de Gr\"obner
de l'annulateur de $f_1^{s_1} \cdots f_p^{s_p}$ dans $\D_{n+p}$,
cf. \cite{bahloulC}.

Pour finir, d\'ecrivons le contenu de l'article.
Dans la premi\`ere section, nous montrons comment, \`a l'aide
de l'op\'erateur d'Euler naturellement associ\'e \`a $\alpha_i$,
faire monter le poids $\rho_{\alpha_i}$ \`a l'aide d'une forme
affine. Dans la section 2, nous d\'ecrivons les
id\'eaux de colongueur finie qui nous servirons dans la suite.
Nous donnons aussi une borne inf\'erieure concernant le poids
qu'il faut pour appartenir \`a ces id\'eaux; ici nous utilisons
la notion de bases standard.
Dans la section 3 se trouve la construction proprement dite
d'un polyn\^ome de Bernstein-Sato. Nous utilisons d'abord
les deux formes affines d\'etermin\'ees dans la section 1 pour faire
monter suffisamment le poids des coefficients des op\'erateurs.
Puis, via une division par les id\'eaux de la section 2, nous faisons
une r\'e\'ecriture de ces op\'erateurs avec une conservation des
poids. Dans la section 4, nous montrons la proposition \ref{propi1}
en utilisant les bases standard. Enfin la section 5 contient la
d\'emonstration de la proposition \ref{propi2}.\\
\\
{\bf Remerciements.} Je tiens \`a
remercier chaleureusement H\'el\`ene Maynadier-Gervais pour les
discussions fructueuses que nous avons eues durant l'\'elaboration de
ce travail. Un grand merci \'egalement \`a Michel Granger pour m'avoir
sugg\'er\'e de tenter le calcul de $\mathcal{H}_f$ et pour avoir relu
le manuscrit. Les r\'esultats expos\'es ici sont issus d'un travail
r\'ealis\'e en grande partie \`a l'universit\'e d'Angers. La section
5 a \'et\'e \'elabor\'ee \`a l'universit\'e de Kobe avec le soutien
financier de la JSPS dans le cadre d'une bourse post-doctorale FY2003.

\section{Processus de mont\'ee des poids}

D\'efinissons les op\'erateurs d'Euler suivants~:
\begin{itemize}
\item $\chi_1=b x_1 \dx{1}+a x_2 \dx{2}$,
\item $\tilde{\chi}_1=b \dx{1} x_1+a \dx{2} x_2=\chi_1+a+b$,
\item $\chi_2=d x_1 \dx{1}+c x_2 \dx{2}$,
\item $\tilde{\chi}_2=d \dx{1} x_1+c \dx{2} x_2=\chi_2+c+d$.
\end{itemize}
Notons alors
\begin{itemize}
\item $f_{11}=\chi_1(f_1)- \rhoun(f_1) f_1$ et $f_{12}=\chi_2(f_1)-
\rhode(f_1) f_1$,
\item $f_{21}=\chi_1(f_2)- \rhoun(f_2) f_2$ et $f_{22}=\chi_2(f_2)-
\rhode(f_2) f_2$.
\end{itemize}
L'affirmation suivante est triviale.

\begin{aff}
\gauche
\begin{itemize}
\item $\rhoun(f_{11}) > \rhoun(f_1)$ et $\rhode(f_{11}) > \rhode(f_1)$,
\item $\rhoun(f_{12}) \ge \rhoun(f_1)$ et $\rhode(f_{12}) >
\rhode(f_1)$,
\item $\rhoun(f_{21}) > \rhoun(f_2)$ et $\rhode(f_{21}) \ge
\rhode(f_2)$,
\item $\rhoun(f_{22}) > \rhoun(f_2)$ et $\rhode(f_{22}) > \rhode(f_2)$.
\end{itemize}
\end{aff}

\begin{nota}
On note $\xi_{-1,0}=f_1^{s_1+1}f_2^{s_2}$, $\xi_{0,-1}=f_1^{s_1}
f_2^{s_2+1}$ et $\xi_{-1,-1}=f_1^{s_1+1}f_2^{s_2+1}$. De plus, si
$(i_1,i_2) \in \N^2$, on note~:
\[\xi_{i_1,i_2}= s_1 \cdots (s_1-i_1+1) \cdot s_2 \cdots (s_2-i_2+1)
f_1^{s_1-i_1} f_2^{s_2-i_2}.\]
\end{nota}

Avant de poursuivre, d\'efinissons le poids d'un op\'erateur
appliqu\'e \`a $\xi_{i_1,i_2}$.

\begin{defin}\label{def:poids}
Soit $i \in \{1,2\}$.
\begin{itemize}
\item
Pour $(i_1,i_2) \in (\N\cup \{-1\})^2$, on pose
$\rho_{\alpha_i}(\xi_{i_1,i_2})=-i_1 \rho_{\alpha_i}(f_1) -i_2
\rho_{\alpha_i}(f_2)$.
\item
Si $u\in \C\{x\}$, on pose $\rho_{\alpha_i}(u \xi_{i_1,i_2})=
\rho_{\alpha_i}(u)+ \rho_{\alpha_i}(\xi_{i_1,i_2})$.
\item
Si $P=P(s) \in \D[s_1,s_2]$, on adopte l'\'ecriture \`a droite~:
$P=\sum_{\beta, k,l} \ddx^\beta s_1^k s_2^l u_{\beta,k,l}$ avec
$\beta\in \N^2$ et $u_{\beta,k,l} \in \C\{x\}$ et on d\'efinit
$\rho_{\alpha_i}(P \cdot \xi_{i_1,i_2})$ comme le minimum des
$\rho_{\alpha_i} (u_{\beta, k,l} \xi_{i_1,i_2})$.
\end{itemize}
\end{defin}

Soient maintenant $(i_1,i_2) \in \N^2$ et $u \in \C\{x\}$ pour lequel
on note $\rho_1=\rhoun(u \xi_{i_1,i_2})$ et $\rho_2=\rhode(u
\xi_{i_1,i_2})$. Posons $u_1=\chi_1(u)- \rhoun(u)u$ et $u_2=\chi_2(u)-
\rhode(u)u$.

\begin{lem}
\begin{eqnarray*}
\tilde{\chi}_1 \cdot u \xi_{i_1,i_2} & =  & \Big( \rhoun(f_1) s_1 +
\rhoun(f_2) s_2 +a+b+\rho_1 \Big) u \xi_{i_1,i_2}\\
&  & + u_1 \xi_{i_1,i_2} + f_{11}u \xi_{i_1+1,i_2} +
f_{21}u \xi_{i_1,i_2+1},
\end{eqnarray*}

\begin{eqnarray*}
\tilde{\chi}_2 \cdot u \xi_{i_1,i_2} & =  & \Big( \rhode(f_1) s_1 +
\rhode(f_2) s_2 +c+d+\rho_2 \Big) u \xi_{i_1,i_2}\\
&  & + u_2 \xi_{i_1,i_2} + f_{12}u \xi_{i_1+1,i_2} +
f_{22}u \xi_{i_1,i_2+1}.
\end{eqnarray*}
\end{lem}
La d\'emonstration consiste en un calcul facile laiss\'e au lecteur.

\begin{rem}
\
\begin{itemize}
\item[$\bullet$] $\rhoun(f_1)s_1+ \rhoun(f_2)s_2+a+b+\rho_1= abs_1+
ads_2+ a+b+ \rho_1$,
\item[$\bullet$] $\rhode(f_1)s_1+ \rhode(f_2)s_2+c+d+\rho_2= ads_1+
cds_2+ c+d+ \rho_2$.
\end{itemize}
\end{rem}

\begin{cor}[de mont\'ee des poids]\label{cor:montee}
Sous les hypoth\`eses du lemme pr\'ec\'edent, on a~:
\begin{enumerate}
\item
$\big( abs_1+ads_2+a+b+\rho_1\big) u\xi_{i_1,i_2} = P_0 \xi_{i_1,i_2}
+ P_1 \xi_{i_1+1,i_2} + P_2 \xi_{i_1,i_2+1},$\\
o\`u $P_0, P_1,P_2 \in \D$ v\'erifient
\[\rhoun(P_0 \xi_{i_1,i_2}), \, \rhoun(P_1 \xi_{i_1+1,i_2}),\,
\rhoun(P_2 \xi_{i_1,i_2+1}) \, > \, \rhoun(u\xi_{i_1,i_2}) \text{ et}\]
\[\rhode(P_0 \xi_{i_1,i_2}), \, \rhode(P_1 \xi_{i_1+1,i_2}),\,
\rhode(P_2 \xi_{i_1,i_2+1}) \, \ge \, \rhode(u\xi_{i_1,i_2}),\]
\item
$\big( ads_1+cds_2+c+d+\rho_2\big) u\xi_{i_1,i_2} = Q_0 \xi_{i_1,i_2}
+ Q_1 \xi_{i_1+1,i_2} + Q_2 \xi_{i_1,i_2+1},$\\
o\`u $Q_0, Q_1,Q_2 \in \D$ v\'erifient
\[\rhode(Q_0 \xi_{i_1,i_2}), \, \rhode(Q_1 \xi_{i_1+1,i_2}),\,
\rhode(Q_2 \xi_{i_1,i_2+1}) \, > \, \rhode(u\xi_{i_1,i_2}) \text{ et}\]
\[\rhoun(Q_0 \xi_{i_1,i_2}), \, \rhoun(Q_1 \xi_{i_1+1,i_2}),\,
\rhoun(Q_2 \xi_{i_1,i_2+1}) \, \ge \, \rhoun(u\xi_{i_1,i_2}).\]
\end{enumerate}
\end{cor}

\section{Les id\'eaux de colongueur f{}inie}

Consid\'erons les id\'eaux suivants de $\C\{x_1,x_2\}$.
\begin{itemize}
\item $I_1=\langle f_1, J \rangle$,
\item $I_2=\langle f_2, J \rangle$,
\item $I=\langle f_1, f_2 \rangle$.
\end{itemize}

Dans cette section, nous d\'eterminons une borne inf\'erieure concernant
le poids d'une s\'erie pour que celle-ci appartienne \`a l'un de
ces id\'eaux. Les r\'esultats
concernant les id\'eaux $I_1$ et $I_2$ sont de H. Maynadier
\cite{maynadierTez, maynadierBull} et ceux qui concernent $I$ s'appuient
sur la notion de bases standard. \`A ce propos nous rappelons
un th\'eor\`eme de division dans le paragraphe 2.

\subsection{Les id\'eaux $I_1$ et $I_2$}
\

Rappelons les hypoth\`eses (5) et (6), \`a savoir $\inu(I_1)=\langle
\inu(f_1), \inu(J) \rangle$ et $\ind(I_2)=\langle \ind(f_2), \ind(J)
\rangle$ et ces id\'eaux sont de colongueur finie. On est alors en
mesure d'appliquer le th\'eor\`eme de division par un id\'eal de
colongueur finie (th. 4.2.1.1 \cite{maynadierBull}) dont on trouvera une
d\'emonstration dans \cite{maynadierTez}, voir aussi \cite{bgmm} et
\cite{briancon}). De plus, la colongueur de $I_1$ (resp. $I_2$) est
finie et \'egale celle de $\langle \inu(f_1), \inu(J) \rangle$
(resp. $\langle \ind(f_2), \ind(J) \rangle$). Des r\'esultats de
H. Maynadier \cite{maynadierTez, maynadierBull}, on retiendra~:

\begin{prop}\label{prop:hlm}
Soit $\mathcal{M}_i$ une cobase monomiale de $I_i$, $i=1,2$. Soit $u$
dans $\C\{x\}$.
\begin{itemize}
\item
Si pour tout $m\in \mathcal{M}_1$, $\rhoun(u)> \degun(m)$,
alors il existe $v,w \in \C\{x\}$ tels que~:
$u=v f_1 + w J$, $\rhoun(u) \le \rhoun(v)+\rhoun(f_1)$ et $\rhoun(u) \le
\rhoun(w)+\rhoun(J)$.
\item
Si pour tout $m\in \mathcal{M}_2$, $\rhode(u)> \degde(m)$, alors il
existe $v,w \in \C\{x\}$ tels que~: $u=v f_2 + w J$, $\rhode(u) \le
\rhode(v)+\rhoun(f_2)$ et $\rhode(u) \le \rhode(w)+\rhoun(J)$.
\end{itemize}
\end{prop}

\begin{rem}\label{rem:hlm}
Dans la th\`ese de H.~Maynadier (\cite{maynadierTez} pages 88-89), il
est expliqu\'e et d\'emontr\'e comment, \'etant donn\'e un id\'eal $J
\subset \C\{x_1,\ldots,x_n\}$ de d\'efinition quasi-homog\`ene pour un
syst\`eme de poids $\alpha \in \N^n$ (par exemple $\inu(I_1)$ ou
$\ind(I_2)$) engendr\'e par des polyn\^omes quasi-homog\`enes
$F_1,\ldots,F_r$ de poids respectifs $p_1,\ldots,p_r$, trouver un
majorant de l'entier $N(J)$ d\'efini comme le plus petit entier $N$
tel que $J$ contienne tous les mon\^omes de poids strictement
sup\'erieur \`a $N$. En particulier lorsque $r=n$ (ce qui est le cas
des id\'eaux $\inu(I_1)$ et $\ind(I_2)$ pour lesquels $r=n=2$),
on a~:
\[N(J) = p_1+\cdots+p_n- |\alpha|.\]
On obtient donc~:
\begin{itemize}
\item[$\bullet$] $N(\inu(I_1))=\rhoun(f_1)+\rhoun(J)-(a+b)$,
\item[$\bullet$] $N(\ind(I_2))=\rhode(f_2)+\rhode(J)-(c+d)$.
\end{itemize}
\end{rem}

\subsection{Pr\'eliminaire sur la division}
\

Dans le paragraphe suivant ainsi que dans la d\'emonstration de
l'affirmation \ref{aff:prop1}, nous allons utiliser la notion
de bases standard dans les s\'eries convergentes $\C\{x\}=
\C\{x_1,\ldots,x_n\}$.
Afin de rendre l'exposition pr\'ecise, nous rappelons le
th\'eor\`eme de division avec uniques reste et quotients
(voir \cite[Th. 1.5.1]{cg}).
Pour les autres r\'esultats et les autres notions, nous renverrons
le lecteur \`a \cite{cg}.

Fixons une forme lin\'eaire $L:\R^n \to \R$ \`a coefficients positifs
ou nuls et soit $<_0$ un bon ordre total fix\'e sur $\N^n$ (par
exemple un ordre lexicographique). On d\'efinit l'ordre $<_L$ sur
$\N^n$ par~:
\[ \beta <_L \beta' \iff L(\beta) > L(\beta') \text{ ou } \big(
\text{\'egalit\'e et } \beta>_0 \beta' \big).\] Pour $f\in
\C\{x\}\smallsetminus\{0\}$ s'\'ecrivant $f=\sum_\beta c_\beta
x^\beta$ (avec $c_\beta \in \C$) on d\'efinit son diagramme de
Newton $\DN(f)=\{\beta \in \N^n|c_\beta \ne 0\}$ puis son
$<_L$-exposant privil\'egi\'e $\exp_{<_L}(f)=\max_{<_L} \DN(f)$.
On note $\cp_{<_L}(f)=c_{\exp_{<_L}(f)}$ son coefficient privil\'egi\'e,
i.e. et $\tp_{<_L}(f)=c_{\exp_{<_L}(f)} x^{\exp_{<_L}(f)}$ son terme
privil\'egi\'e. Finalement, on note $\rho_L(f)$ son poids d\'efini comme
le minimum de l'image $L(\DN(f))$.
Signalons que $\rho_L(f)= L(\exp_{<_L}(f))$.

Soient maintenant $f_1,\ldots,f_r$ dans
$\C\{x\}\smallsetminus\{0\}$. Consid\'erons la partition
$\N^n=\Delta_1 \cup \cdots \cup \Delta_r \cup \bar{\Delta}$ associ\'ee
\`a $f_1, \ldots, f_r$~:
\begin{itemize}
\item $\Delta_1=\exp_{<_L}(f_1) +\N^n$,
\item $\Delta_j=(\exp_{<_L}(f_j) +\N^n) \smallsetminus (\Delta_1 \cup
\cdots \cup \Delta_{j-1})$, pour $j=2,\ldots, r$,
\item $\bar{\Delta}=\N^n \smallsetminus (\Delta_1 \cup \cdots \cup
\Delta_r)$.
\end{itemize}

\begin{theo}[de division]
Pour tout $f\in \C\{x\}$, il existe un unique \'el\'ement
$(q_1,\ldots,q_r,R)$ de $(\C\{x\})^{r+1}$ tel que~:
\begin{enumerate}
\item $f=q_1 f_1 +\cdots + q_r f_r + R$,
\item $\DN(q_j) +\exp_{<_L}(f_j) \subset \Delta_j$ pour tout $j$ tel
que $q_j \ne 0$,
\item $\DN(R) \subset \bar{\Delta}$ si $R \ne 0$.
\end{enumerate}
L'\'el\'ement $R$ est appel\'e le reste de la division.
\end{theo}

Rappelons les grandes lignes de la d\'emonstration du th\'eor\`eme
de division. Cela nous sera utile dans la suite.
\begin{description}
\item[(i)] Posons $(f^0, q_1^0, \ldots, q_r^0, R^0)=(f, 0, \ldots,0,
0)$.
\end{description}
Pour $i\ge 0$,
\begin{description}
\item[(ii)] si $f^i=0$ alors on pose $(f^{i+1}, q_1^{i+1}, \ldots, q_r^{i+1}, R^{i+1})=(f^i, q_1^i, \ldots, q_r^i, R^i)$.
\item[(iii)] si $\exp_{<_L}(f^i) \in \bar{\Delta}$ alors
\[(f^{i+1}, q_1^{i+1}, \ldots, q_r^{i+1}, R^{i+1})=
(f^i-\tp_{<_L}(f^i), q_1^{i+1}, \ldots, q_r^{i+1}, R^i+
\tp_{<_L}(f^i)),\]
\item[(iv)] sinon, soit $j=\min\{k \in \{1,\ldots,r\} ,\,
\exp_{<_L}(f^i) \in \Delta_k\}$ et

$\dps f^{i+1}=f^i- \frac{\cp_{<_L}(f^i)}{\cp_{<_L}(f_j)} \cdot
x^{\exp_{<_L}(f^i)-\exp_{<_L}(f_j)} \cdot f_j$,\\

$\dps q_j^{i+1}=q_j^i + \frac{\cp_{<_L}(f^i)}{\cp_{<_L}(f_j)} \cdot
x^{\exp_{<_L}(f^i)-\exp_{<_L}(f_j)}$,\\

$\dps q_l^{i+1}=q_l^i$ pour $l\ne j$, et $\dps R^{i+1}=R^i$.
\end{description}
Cette construction donne lieu a $r+2$ suites $f^i, q_1^i, \ldots,
q_r^i, R^i$ qui v\'erifient $f=f^i+\sum_{j=1}^r q_j^i f_j +R^i$ et
dont on montre qu'elles convergent pour la topologie
$(x_1,\ldots,x_n)$-adique (en particulier les $f^i$ tendent vers
$0$). On obtient donc une division dans $\C[[x]]$. La fin de la preuve
consiste \`a montrer que $R$ ainsi que les $q_j$ construits sont bien
dans $\C\{x\}$. On appelle \emph{division \'el\'ementaire} l'une des
\'etapes (ii), (iii), (iv) dans la construction ci-dessus.\\
Pour finir, voici une remarque utilis\'ee dans la d\'emonstration du
lemme \ref{lem:bsI} et de l'affirmation \ref{aff:prop1}.

\begin{rem}\label{rem:utile}
Dans la construction ci-dessus, si pour tout $i$, $\exp_{<_L}(f^i)$
appartient \`a $\Delta_1\cup \cdots\cup \Delta_r$ alors le reste $R$
est nul (en effet chaque $R^i$ est nul).
\end{rem}

\subsection{L'id\'eal $I$}
\

Dans ce paragraphe, nous montrons que l'id\'eal $I$ est de colongueur
finie; ceci gr\^ace \`a une division relativement aux deux syst\`emes de
poids $\alpha_1$ et $\alpha_2$. Pour cela, nous allons utiliser
la notion de bases standard (voir \cite{cg}).

Soit $L_1$ (resp. $L_2$) la forme lin\'eaire $L_1(i,j)=bi+aj$
(resp. $L_2(i,j)=di+cj$). Soient $<_1$ et $<_2$ les ordres sur $\N^2$
d\'efinis comme suit~:
\[ (i,j) <_1 (i',j') \iff L_1(i,j) > L_1(i',j') \text{ ou } \big(
\text{\'egalit\'e et } j>j' \big),\]
\[ (i,j) <_2 (i',j') \iff L_2(i,j) > L_2(i',j') \text{ ou } \big(
\text{\'egalit\'e et } i>i' \big).\] Pour $u \in \C\{x\}$, on note
$\exp_{<_1}(u)=\max_{<_1}(\DN(u))$ (de m\^eme pour $<_2$).  Nous avons
les propri\'et\'es suivantes~:

\begin{propri}\label{propri}
\
\begin{description}
\item[(i)] $L_1(\exp_{<_1}(u))=\rhoun(u)$ et
  $L_2(\exp_{<_2}(u))=\rhode(u)$
\item[(ii)] $\exp_{<_1}(f_1)=\exp_{<_2}(f_1)=(a,0)$
\item[(iii)] $\exp_{<_1}(f_2)=\exp_{<_2}(f_2)=(0,d)$
\end{description}
\end{propri}

Dans ce qui suit, nous montrons que $f_1$ et $f_2$ forment une
$<_1$-base standard de $I$. Gr\^ace aux propri\'et\'es {\bf (ii)} et
{\bf (iii)}, cela entraine qu'ils forment aussi une $<_2$-base
standard.

\begin{lem}\label{lem:bsI}
L'ensemble $\{f_1,f_2\}$ est une $<_1$-base standard de $I$.
\end{lem}
\begin{proof}
Pour cela, nous allons consid\'erer la division de la $S$-fonction
(cf. \cite[Def. 1.6.1]{cg}),
$S=x_2^d f_1 -x_1^a f_2$, par $f_1,f_2$ relativement \`a l'ordre
$<_1$~: $S=q_1 f_1 + q_2 f_2 +R$, le but \'etant de montrer que le
reste $R$ est nul (il s'agit du crit\`ere de Buchberger d\'emontr\'e
\`a l'origine dans les anneaux de polyn\^omes \cite{buchberger}).
Dans $\N^2$, consid\'erons le secteur \'epoint\'e
\[\Gamma=\{(i,j)/\, L_1(i,j)\ge L_1(a,d) \text{ et }
L_2(i,j)\ge L_2(a,d)\} \smallsetminus \{(a,d)\}.\] On remarque que
$\DN(S) \subset \Gamma$ et surtout que
\[\Gamma \subset (\exp_{<_1}(f_1) +\N^2) \cup (\exp_{<_1}(f_2)
+\N^2).\]  Consid\'erons maintenant la division de $S$ par $(f_1,f_2)$
relativement \`a $<_1$. Formellement, elle consiste en une suite
infinie de divisions \'el\'ementaires.
La premi\`ere d'entre elles consiste \`a
diviser le mon\^ome privil\'egi\'e $m=\exp_{<_1}(S)=x_1^kx_2^l$ de $S$
par $f_1$ (ou $f_2$). En terme de diagramme de Newton, cette division
se traduit par le fait de remplacer $(k,l)$ par l'ensemble
\[\big( \DN(f_1) + (k-a,l)\big) \smallsetminus \{(k,l)\}\]
(ou $\big( \DN(f_2) + (k,l-d)\big) \smallsetminus \{(k,l)\}$) qui est
inclus dans $\Gamma$.  On voit alors que tout au long du processus de
division, toutes les divisions \'el\'ementaires se d\'eroulent dans le
secteur $\Gamma$ (i.e. avec les notations employ\'ees dans le
paragraphe pr\'ec\'edent, pour tout $i \in \N$, $\DN(f^i)$ est inclus
dans $\Gamma$). En cons\'equence, le reste de la division de $S$ par
$(f_1,f_2)$ est n\'ecessairement nul (cf. remarque \ref{rem:utile}).
En utilisant le crit\`ere de Buchberger \cite[Prop. 1.6.2]{cg},
on peut conclure que $f_1,f_2$ forment une $<_1$-base standard de $I$.
\end{proof}

\begin{cor}
L'id\'eal $I$ est de colongueur finie et $\mathcal{M}=\{x_1^ix_2^j/ \,
0\le i \le a-1,\, 0\le j \le d-1\}$ en est une cobase.
\end{cor}

\begin{cor}\label{cor:divI}
Soit $u \in \C\{x\}$ tel que $\rhoun(u)> \degun(x_1^{a-1} x_2^{d-1})$
ou $\rhode(u)> \degde(x_1^{a-1} x_2^{d-1})$ alors il existe $u_1,u_2
\in \C\{x\}$ tels que~:
\begin{itemize}
\item[$\bullet$] $u=u_1 f_1 +u_2 f_2$,
\item[$\bullet$] $\rhoun(u) \le \rhoun(u_1) +\rhoun(f_1)$,\\
$\rhoun(u) \le \rhoun(u_2) +\rhoun(f_2)$,
\item[$\bullet$] $\rhode(u) \le \rhode(u_1) +\rhode(f_1)$,\\
$\rhode(u) \le \rhode(u_2) +\rhode(f_2)$.
\end{itemize}
\end{cor}
\begin{proof}
Divisons $u$ par $f_1,f_2$ relativement \`a $<_1$. On obtient~: $u=u_1
f_1 +u_2 f_2 + R$, avec~:
\begin{eqnarray*}
\rhoun(u) & \le & \rhoun(u_1) +\rhoun(f_1),\\
\rhoun(u) & \le & \rhoun(u_2) +\rhoun(f_2)\\
\text{et } \rhoun(u) & \le & \rhoun(R).
\end{eqnarray*}
Or, d'apr\`es les propri\'et\'es \ref{propri} {\bf (ii)} et {\bf (iii)},
la division pr\'ec\'edente est aussi une division relativement \`a
$<_2$, on a donc~:
\begin{eqnarray*}
\rhode(u) & \le & \rhode(u_1) +\rhode(f_1),\\
\rhode(u) & \le & \rhode(u_2) +\rhode(f_2)\\
\text{et } \rhode(u) & \le & \rhode(R).
\end{eqnarray*}
Comme $R$ est une combinaison lin\'eaire de mon\^omes de
$\mathcal{M}$, il est n\'ecessairement nul.
\end{proof}

\section{Construction d'un polyn\^ome de Bernstein-Sato}

Ce qui suit contient la construction et la
formule explicite d'un polyn\^ome de Bernstein-Sato associ\'e \`a
$(f_1,f_2)$. L'id\'ee g\'en\'erale est la suivante : en appliquant
les formes affines obtenues en section 1, nous faisons une mont\'ee
suffisante des poids d'op\'erateurs agissant sur $\xi_{i_1, i_2}$.
Ici suffisante signifie qu'on veut que les coefficients de ces
op\'erateurs appartiennent aux id\'eaux $I$, $I_1$ et $I_2$.
Ensuite via une division par ces id\'eaux, on montre comment passer
de $\xi_{i_1,i_2}$ \`a $\xi_{i_1-1, i_2}$ et $\xi_{i_1, i_2-1}$
ceci en conservant le(s) poids des op\'erateurs, le but \'etant
de ``redescendre'' jusqu'\`a $\xi_{-1,-1}=f_1^{s_1+1}f_2^{s_2+1}$.

Pour commencer, faisons quelques remarques sur $J$ le d\'eterminant
jacobien de $(f_1,f_2)$.
\begin{rem}\label{rem:surJ}
\gauche
\begin{itemize}
\item $(a-1,d-1)$ appartient \`a $\DN(J)$,
\item $\dps \rhoun(J)=\degun(x_1^{a-1}x_2^{d-1})=b(a-1)+a(d-1)$,
\item $\dps \rhode(J)=\degde(x_1^{a-1}x_2^{d-1})=d(a-1)+c(d-1)$,
\item $\dps \inu(J)=\inu(\dxsur{f_2}{2}) \cdot \inu(\dxsur{f_1}{1})=d
  x_2^{d-1} \cdot \inu(\dxsur{f_1}{1})$,
\item $\dps \ind(J)=\ind(\dxsur{f_1}{1}) \cdot \ind(\dxsur{f_2}{2})=a
  x_1^{a-1} \cdot \ind(\dxsur{f_2}{2})$.
\end{itemize}
\end{rem}
La preuve est facile et laiss\'ee au lecteur.\\

Notons $N(I_1)=\rhoun(f_1)+\rhoun(J)-(a+b)$ et $N_1(I)=b(a-1)+
a(d-1)$. Nous avons vu que si un mon\^ome $m$ a un $\alpha_1$-degr\'e
strictement sup\'erieur \`a $N(I_1)$ alors il appartient \`a $I_1$
(voir la proposition \ref{prop:hlm} et la remarque \ref{rem:hlm}). De
m\^eme si son $\alpha_1$-degr\'e est strictement sup\'erieur \`a $N_1(I)$,
il appartient \`a $I$ (voir le corollaire \ref{cor:divI}). Une
question naturelle est celle de comparer $N(I_1)$ et $N_1(I)$ (et de
m\^eme pour $N(I_2)= \rhode(f_2)+\rhode(J)-(c+d)$ et
$N_2(I)=d(a-1)+c(d-1)$).

\begin{lem}\label{lem:comparaison}
Si $a\ge 2$ (resp. $d\ge 2$) alors $N(I_1) \ge N_1(I)$ (resp. $N(I_2)
\ge N_2(I)$).
\end{lem}
La preuve de ce lemme est laiss\'ee au lecteur.

\begin{lem}\label{lem:lisse}
Si $a=1$ et $d=1$ alors l'id\'eal de Bernstein-Sato de $(f_1,f_2)$ est
\'egal \`a $\langle (s_1+1)(s_2+1) \rangle$.
\end{lem}

\begin{proof}
Par la remarque pr\'ec\'edente \ref{rem:surJ}, $(a-1,d-1) \in \DN(J)$
ainsi $J$ est une unit\'e de $\C\{x\}$. Il en r\'esulte alors que
$(s_1+1) (s_2+1)$ est un polyn\^ome de Bernstein-Sato.  Il est bien
connu d'autre part que puisque $(f_1,f_2)(0)=(0,0)$ et que $f_1,f_2$
sont en intersection compl\`ete (en effet $I$ est de colongueur finie
d'o\`u $V(f_1,f_2)=V(I)=(0)$), $\B$ est inclus dans $\langle (s_1+1)
(s_2+1) \rangle$ (voir le lemme 1.2 dans \cite{maynadierBull}).
\end{proof}

Ainsi, dans la suite nous supposerons toujours~: $a\ge 2$ ou $d\ge 2$.

\subsection{\'Etape 1~: Mont\'ee suffisante des poids $\rho_1$ et $\rho_2$}

Consid\'erons les deux ensembles suivants~:

\[W_1=\{\degun(m)/ \, m \text{ mon\^ome et } \degun(m)\le N(I_1)+
\rhoun(f_2)\},\]
\[W_2=\{\degde(m)/ \, m \text{ mon\^ome et } \degde(m)\le N(I_2)+
\rhode(f_1)\}.\] Ensuite, posons
\[\tilde{b}(s_1,s_2)= \prod_{\rho_1 \in W_1} (abs_1 +ads_2+ a+b+
\rho_1) \cdot \prod_{\rho_2 \in W_2} (ads_1 +cds_2+ c+d+ \rho_2),\] et
enfin $b(s_1,s_2)=(s_1+1)(s_2+1) \tilde{b}(s_1,s_2)$. Le but est de
montrer que $b(s_1,s_2)$ est un polyn\^ome de Bernstein-Sato de
$(f_1,f_2)$.

\begin{lem}\label{lem:etape1}
Il existe des entiers $M_1$ et $M_2$ et des op\'erateurs $P_{i_1,i_2}
\in \D$ pour $i_1=0,\ldots, M_1$ et $i_2=0,\ldots,M_2$ tels que~:
\[\tilde{b}(s_1,s_2) f_1^{s_1}f_2^{s_2} = \sum_{\substack{0\le i_1 \le M_1\\
0 \le i_2 \le M_2}} P_{i_1,i_2} \cdot \xi_{i_1,i_2}\] et
pour tout $i_1,i_2$,
\[\rhoun(P_{i_1,i_2} \cdot \xi_{i_1,i_2}) > N(I_1)+\rhoun(f_2) \quad
\text{et} \quad \rhode(P_{i_1,i_2} \cdot \xi_{i_1,i_2}) >
N(I_2)+\rhode(f_1).\]
\end{lem}
Ce lemme d\'ecoule directement du corollaire \ref{cor:montee}.

\subsection{\'Etape 2~: Division par $I$ et r\'e\'ecriture}

Dans le lemme pr\'ec\'edent, pour tout $i_1$ et tout $i_2$,
$P_{i_1,i_2}$ est une somme finie de $\dx{1}^{\beta_1} \dx{2}^{\beta_2}
u$ avec des $u \in \C\{x\}$ qui v\'erifient~:

$\rhoun(u) > N(I_1)+\rhoun(f_2) +i_1 \rhoun(f_1) +i_2 \rhoun(f_2)$ et

$\rhode(u) > N(I_2)+\rhode(f_1) +i_1 \rhode(f_1) +i_2 \rhode(f_2)$.\\
Ceci nous dit en particulier que chaque $u$ appartient
\`a $I$ (cf. lemme \ref{lem:comparaison} et corollaire
\ref{cor:divI}).  Divisons un tel $u$ par $(f_1,f_2)$ relativement \`a
l'ordre $<_1$ (et aussi $<_2$, voir la propri\'et\'e \ref{propri}). On
obtient $u=u_1 f_1 + u_2 f_2$ avec les conditions suivantes~:
\[\rhoun(u) \le \rhoun(u_1) + \rhoun(f_1) \quad \text{et} \quad
\rhoun(u) \le \rhoun(u_2) + \rhoun(f_2),\]
\[\rhode(u) \le \rhode(u_1) + \rhode(f_1) \quad \text{et} \quad
\rhode(u) \le \rhode(u_2) + \rhode(f_2).\] En cons\'equence~:

\begin{lem}[R\'e\'ecriture et conservation des poids]
Avec les notations du lemme pr\'ec\'edent \ref{lem:etape1}, soit
$(i_1,i_2) \in \{0,\ldots, M_1 \} \times \{0,\ldots,M_2 \}$. Il existe
$P^1, P^2 \in \D$ tels que~:
\[P_{i_1,i_2} \cdot \xi_{i_1,i_2} = (s_1-i_1+1)P^1 \cdot
\xi_{i_1-1,i_2} + (s_2-i_2+1)P^2 \cdot \xi_{i_1,i_2-1}\]  et
\[\rhoun(P^1 \cdot \xi_{i_1-1,i_2}) > N(I_1)+\rhoun(f_2) \quad
\text{et} \quad \rhode(P^1 \cdot \xi_{i_1-1,i_2}) > N(I_2)+
\rhode(f_1),\] et de m\^eme pour $P^2 \cdot \xi_{i_1,i_2-1}$.
\end{lem}
En it\'erant ce processus de division et r\'e\'ecriture, on obtient
ceci~:
\begin{cor}\label{cor:etape2}
Il existe des entiers $M'_1,M'_2$ et des op\'erateurs $Q_{i_1},
R_{i_2} \in \D[s_1,s_2]$ avec $0\le i_1\le M'_1 \le M_1$ et $0\le i_2
\le M'_2 \le M_2$ tels que
\[\tilde{b}(s_1,s_2) f_1^{s_1}f_2^{s_2} = \sum_{i_1=0}^{M'_1} Q_{i_1}
\cdot \xi_{i_1,-1} + \sum_{i_2=0}^{M'_2} R_{i_2} \cdot \xi_{-1,i_2}\]
et pour tout $i_1$ et tout $i_2$,
\[\rhoun(Q_{i_1} \cdot \xi_{i_1,-1}) > N(I_1) +\rhoun(f_2) \quad
\text{et} \quad \rhode(R_{i_2}  \cdot \xi_{-1,i_2}) > N(I_2)
+\rhode(f_1).\]
\end{cor}

En adoptant l'\'ecriture \`a droite, chaque $Q_{i_1}$ (resp. $R_{i_2}$)
peut s'\'ecrire comme une somme finie de $\ddx^\beta s_1^k s_2^lu$,
$u\in \C\{x\}$ (resp. $\ddx^\beta s_1^k s_2^l v$, $v\in \C\{x\}$) de
telle sorte que~:
\begin{enumerate}
\item
$\rhoun(u) > N(I_1) +i_1 \rhoun(f_1)$,
\item
$\rhode(v) > N(I_2) +i_2 \rhode(f_2)$.
\end{enumerate}

L'\'etape suivante va consister \`a diviser par $I_1$ les $u$
satisfaisant (1) et par $I_2$ les $v$ satisfaisant (2).

\subsection{\'Etape 3~: Division par $I_1$ (resp. $I_2$) et r\'e\'ecriture}

On se contentera de traiter le premier cas (les calculs sont
similaires dans le second).

Soit donc $u \in \C\{x\}$ tel que $\rhoun(u) > N(I_1) +i_1
\rhoun(f_1)$, $i_1$ \'etant fix\'e sup\'erieur ou \'egal \`a $1$ (si
$i_1=0$, on ne fait rien).

Par hypoth\`ese sur son $\alpha_1$-poids, $u$ appartient \`a $I_1$ et
par division, on peut \'ecrire~: $u=v f_1 + w J$ avec
\[ \rhoun(u) \le \rhoun(v) + \rhoun(f_1) \quad \text{et} \quad
\rhoun(u) \le \rhoun(w) + \rhoun(J).\] Ainsi,
\[u\,\xi_{i_1,-1} = (s_1-i_1+1) v \xi_{i_1-1,-1} + w J \xi_{i_1,-1}\]
o\`u l'on a $\rhoun(v) > N(I_1)+ (i_1-1)\rhoun(f_1)$.\\
Nous allons maintenant nous occuper du terme $w J \xi_{i_1,-1}$.

\begin{aff}
\[J \, \xi_{i_1,-1}=\Big( \dx{1} \dxsur{f_2}{2} - \dx{2}
\dxsur{f_2}{1}\Big) \cdot \xi_{i_1-1,-1}.\]
\end{aff}
La d\'emonstration consiste en un simple calcul laiss\'e au lecteur.

En multipliant cette \'egalit\'e par $w$ et en faisant commuter ce
dernier avec les $\dx{i}$, on obtient~:
\[wJ \xi_{i_1,-1}= \Big( \dx{1} \underbrace{w \dxsur{f_2}{2}}_{E_1}
-\dx{2} \underbrace{w \dxsur{f_2}{1}}_{E_2} + \underbrace{\dxsur{w}{1}
\dxsur{f_2}{2}}_{E_3} - \underbrace{\dxsur{w}{2}
\dxsur{f_2}{1}}_{E_4} \Big)\cdot \xi_{i_1-1,-1}.\]

La suite consiste \`a contr\^oler le $\alpha_1$-poids des termes $E_1,
\ldots,E_4$ ainsi obtenus, le but \'etant d'avoir~: $\rhoun(E_i) >
N(I_1) +(i_1-1) \rhoun(f_1)$.

\begin{aff}
Pour $i=1,\ldots,4$, $\rhoun(E_i) > N(I_1) +(i_1-1) \rhoun(f_1)$.
\end{aff}
\begin{proof}
En premier lieu, remarquons qu'il suffit d'avoir les in\'egalit\'es
pour $E_3$ et $E_4$ (en effet, $\rhoun(E_1) \ge \rhoun(E_3)$ et
$\rhoun(E_2) \ge \rhoun(E_4)$).  Montrons l'in\'egalit\'e concernant
$E_3$. Pour cela, il suffit d'obtenir $\rhoun(E_3) \ge \rhoun(u)
-\rhoun(f_1)$. Nous avons
\begin{eqnarray*}
\rhoun(E_3) & \ge & \rhoun(w)- \degun(x_1)+ \rhoun(f_2)- \degun(x_2)
\qquad (\star)\\
& \ge & \rhoun(u) -\rhoun(J) -\degun(x_1)+\rhoun(f_2)-\degun(x_2)\\
& = & \rhoun(u)- b(a-1) -a(d-1) -b +ad -a\\
& = & \rhoun(u) -ab\\
& = & \rhoun(u) -\rhoun(f_1).
\end{eqnarray*}
L'in\'egalit\'e qui concerne $E_3$ est d\'emontr\'ee. Pour $E_4$, nous
avons~:
\[\rhoun(E_4) \ge \rhoun(w) -\degun(x_2) +\rhoun(f_2) -\degun(x_1),\]
et l'on retombe sur $(\star)$.
\end{proof}

\subsection{\'Etape f\mbox{}inale~: le bilan}

Notons $\mathbb{D}$ l'alg\`ebre des op\'erateurs diff\'erentiels \`a
coefficients constants. Alors, comme cons\'equence des \'etapes
pr\'ec\'edentes, nous avons~:

\[\tilde{b}(s_1,s_2) f_1^{s_1} f_2^{s_2} \in \sum_u
\mathbb{D}[s_1,s_2] u \, \xi_{0,-1} + \sum_v \mathbb{D}[s_1,s_2] v \,
\xi_{-1,0}\] avec pour chaque $u$, $\rhoun(u) > N(I_1)$ et pour chaque
$v$, $\rhode(v) > N(I_2)$.

Autrement dit chaque $u$ appartient \`a $I_1$ et chaque $v$ \`a $I_2$.
Prenons le cas d'un $u$ et \'ecrivons $u=u_1 f_1 +u_2 J$. Par suite,
\[(s_1+1) u\, \xi_{0,-1} = (s_1+1) u_1 \, \xi_{-1,-1} + u_2 \Big(
\dx{1} \dxsur{f_2}{2} - \dx{2} \dxsur{f_2}{1} \Big) \cdot
\xi_{-1,-1}\] ou encore
\[(s_1+1) u \, \xi_{0,-1} \in \D[s_1,s_2] f_1^{s_1+1} f_2^{s_2+1}.\]
De mani\`ere similaire on obtient~:
\[(s_2+1) v \, \xi_{-1,0} \in \D[s_1,s_2] f_1^{s_1+1} f_2^{s_2+1}.\]
En conclusion~:
\[b(s_1,s_2) f_1^{s_1} f_2^{s_2} \in \D[s_1,s_2] f_1^{s_1+1}
f_2^{s_2+1}\] et le th\'eor\`eme \ref{theoi} est d\'emontr\'e.

\section{D\'emonstration de la proposition \ref{propi1}}

Nous devons montrer que $(f_1,f_2)$ satisfait aux hypoth\`eses
(1),\ldots,(6). Les quatre premi\`eres \'etant clairement
v\'erifi\'ees, il nous suffit de le faire pour (5) et (6). Enfin, par
sym\'etrie, il suffit de montrer (5), \`a savoir~:

\begin{aff}\label{aff:prop1}
\
\begin{itemize}
\item[(i)] $I_1$ est de colongueur finie.
\item[(ii)] $\inu(f_1)$ et $\inu(J)$ engendrent $\inu(I_1)$.
\end{itemize}
\end{aff}
L'hypoth\`ese (5) d\'ecoule directement de
\cite[Th. 4.2.1.1]{maynadierBull} (cf. \cite[part. 2]{maynadierTez} pour
les d\'etails).

\begin{proof}
Nous allons d\'emontrer les deux points en m\^eme temps en
d\'eterminant une $<_1$-base standard de $I_1$. Pour cela, nous allons
suivre l'algorithme de Buchberger \cite[1.6]{cg}. Aussi, nous
adopterons la notation suivante~: si $u_1, \ldots,u_m$ sont dans
$\C\{x\}$, nous noterons
\[E(u_1,\ldots,u_m)=\bigcup_{i=1}^m \Big( \exp_{<_1}(u_i) +\N^2
\Big).\] Il est facile de voir que $\exp_{<_1}(f_1)=(a,0)$. De plus,
gr\^ace \`a la remarque \ref{rem:surJ}, $\inu(J)=ad x_1^{a-1}
x_2^{d-1}$ ce qui entraine $\exp_{<_1}(J)=(a-1,d-1)$. Quitte \`a
diviser $J$ par $ad$, on peut supposer que $J=x_1^{a-1} x_2^{d-1}+h$
avec $\rhoun(h) > \rhoun(J)$. Maintenant, consid\'erons la
$S$-fonction de $f_1$ et $J$~:
\[K=S(f_1,J)=x_2^{d-1} f_1 - x_1 J.\]
Nous avons $\inu(x_2^{d-1} f_1)=x_1^a x_2^{d-1}+ x_2^{b+d-1}$ et
$\inu(x_1J)= x_1^a x_2^{d-1}$. Ainsi, $\inu(K)=x_2^{b+d-1}$ ce qui
entraine $\exp_{<_1}(K)=(0,b+d-1)$.\\
Notons qu'\`a ce stade, le point (i) de l'affirmation est acquis
(cf \cite[Cor. 1.5.3]{cg}).\\
Poursuivons en d\'emontrant que
$f_1$, $J$ et $K$ forment une $<_1$-base standard de $I_1$ (remarquons
que $K$ n'est divisible ni par $f_1$, ni par $J$). Nous allons montrer
que le reste de la division de $S(f_1,K)$ (resp. de $S(J,K)$) par
$(f_1,J,K)$ est nul.

Consid\'erons la $S$-fonction $S_1=S(f_1,K)=x_2^{b+d-1} f_1- x_1^a
K$. On constate que $\inu(S_1)=x_2^{2b+d-1}$. Par cons\'equent,
$\DN(S_1)$ est inclus dans le demi-espace ferm\'e~:
\[\Gamma_1=\{(k,l)\in \N^2/\, L_1(k,l) \ge L_1(0,2b+d-1)\}.\]
On voit ais\'ement que $\Gamma_1 \subset E(f_1,J,K)$. Maintenant,
consid\'erons la division de $S_1$ par $(f_1,J,K)$. On voit que toutes
les divisions \'el\'ementaires qui la constituent se font avec des
diagrammes de Newton qui ne sortent pas de $\Gamma_1$ (donc de
$E(f_1,J,K)$). Par cons\'equent, le reste de cette division est nul
(voir la remarque \ref{rem:utile}).

En ce qui concerne $S_2=S(J,K)= x_1^{a-1}K- x_2^b J$. L'argument est
similaire et s'appuie sur le fait que $\DN(S_2)$ est inclus dans le
demi-espace ouvert $\Gamma_2=\{(k,l)\in \N^2/\, L_1(k,l) > L_1(a-1,
b+d-1)\}$, lui-m\^eme inclus dans $E(f_1,J,K)$.

Comme bilan, l'ensemble form\'e de $f_1$, $J$ et $K$ est une
$<_1$-base standard de $I_1$ \cite[Prop. 1.6.2]{cg}.
Une cons\'equence directe est que
l'id\'eal $\inu(I_1)$ est engendr\'e par $\inu(f_1)$, $\inu(J)$ et
$\inu(K)$. Or $\inu(K)=x_2^{b+d-1}$ et $x_2^{d-1} \inu(f_1) - x_1
\inu(J)= x_2^{b+d-1}$. Ainsi $\inu(K) \in \langle \inu(f_1), \inu(J)
\rangle$ ce qui ach\`eve la preuve du point (ii) et celle de
l'affirmation.
\end{proof}

\section{D\'emonstration de la proposition \ref{propi2}}

Rappelons que nous supposons toujours $a\ge 2$ ou $d\ge 2$.\\
Pour la d\'emonstration de la prop. \ref{propi2}, il suffit de montrer les
inclusions $V(s_1 s_2) \subset \mathcal{H}_f$  et $V((ab s_1 + ad s_2)
(ad s_1+ cd s_2)) \subset \mathcal{H}_f$. Par sym\'etrie, nous pouvons
nous contenter de~:
\begin{aff}
$V(s_1s_2) \subset \mathcal{H}_f$ et $V(ab s_1 + ad s_2) \subset
\mathcal{H}_f$.
\end{aff}
\begin{proof}
Commen\c{c}ons par l'inclusion $V(s_1s_2) \subset \mathcal{H}_f$.\\
Comme nous l'avons d\'ej\`a vu, puisque $f_1$ et $f_2$ s'annulent en $0$
et sont en intersection compl\`ete, tout polyn\^ome de Bernstein est inclus
dans $\langle (s_1+1)(s_2+1) \rangle$. Par cons\'equent, $\fin(\B)$ est
multiple de $\langle s_1s_2 \rangle$. L'inclusion recherch\'ee
d\'ecoule alors du th\'eor\`eme BMM rappel\'e dans l'introduction.\\
Montrons l'inclusion $V(ab s_1 + ad s_2) \subset \mathcal{H}_f$.\\
Via l'identification $T^*\C^2 \times \C^2= \C^6$, $A$ est \'egal \`a
l'ensemble des
\[\big( x_1,x_2, \underbrace{\lambda_1 a x_1^{a-1}+ \lambda_2
c x_1^{c-1}}_{\xi_1}, \underbrace{\lambda_1 b x_2^{b-1}+ \lambda_2 d
x_2^{d-1}}_{\xi_2}, \underbrace{\lambda_1 (x_1^a+ x_2^b)}_{s_1},
\underbrace{\lambda_2 (x_1^c+ x_2^d)}_{s_2}\big).\]
Soit $(s_1,s_2)$ tel que $abs_1 +ads_2=0$. Nous allons traiter deux
cas s\'epar\'ement~:\\
Cas 1 : $bc-ad+ d-b \ge 0$.\\
Dans ce cas, on d\'efinit les suites suivantes~:
\begin{itemize}
\item
$x_1(n)$ est une suite de r\'eels positifs non nuls qui tendent vers $0$.
\item
Pour tout $n$, soit $x_2(n)\in \R_{>0}$ tel que
$x_2(n)^b=x_1(n)^{a-1}$.
\item
Pour tout $n$, $\lambda_1(n)= s_1 x_1(n)^{1-a}$.
\item
Pour tout $n$, $\lambda_2(n)=-\frac{b}{d} s_1 x_2(n)^{-d}$.
\end{itemize}
Notons $\xi_1(n)=\lambda_1(n) a x_1(n)^{a-1}+ \lambda_2(n) c
x_1(n)^{c-1}$. De la m\^eme mani\`ere on d\'efinit $\xi_2(n)$, $s_1(n)$ et
$s_2(n)$.
Nous constatons alors que~:
\begin{enumerate}
\item $\xi_1(n)=a s_1 -\frac{b}{d} s_1 x_2(n)^{-d} x_1(n)^{c-1}$,
\item $\xi_2(n)=0$,
\item $s_1(n)= s_1 x_1(n)+ s_1$,
\item $s_2(n)= -\frac{b}{d} s_1 (1 + x_2(n)^{-d} x_1(n)^c)$.
\end{enumerate}
En \'elevant $x_2(n)^{-d} x_1(n)^{c-1}$ \`a la puissance $b$ et en
utilisant la relation $x_2(n)^b=x_1(n)^{a-1}$, on constate que
$(x_2(n)^{-d} x_1(n)^{c-1})^b= x_1(n)^{bc-ad +d-b}$. Ainsi la suite
$x_2(n)^{-d} x_1(n)^{c-1}$ tend vers $0$ si $bc-ad +d-b>0$ et est
constante si $bc-ad +d-b=0$ (c'est sous cette hypoth\`ese
qu'intervient le fait que les suites $x_1(n)$ et $x_2(n)$ sont dans
$\R_{>0}$).  On constate alors que les quatre suites $\xi_1(n),\xi_2(n),
s_1(n),s_2(n)$ sont convergentes et que $(s_1(n), s_2(n))$ converge
vers $(s_1,s_2)$ qui appartient donc bien \`a $\mathcal{H}_f$.\\
Cas 2 : $bc-ad+d-b <0$.\\
Ce cas se traite de la m\^eme fa\c{c}on avec les donn\'ees suivantes (les
d\'etails sont laiss\'es au lecteur).
\begin{itemize}
\item
$x_1(n)$ est une suite de complexes non nuls qui tendent vers $0$.
\item
Pour tout $n$, soit $x_2(n)\in \C$ tel que $x_2(n)^d=x_1(n)^{c-1}$.
\item
Pour tout $n$, $\lambda_2(n)= s_2 x_2(n)^{-d}$.
\item
Pour tout $n$, $\lambda_1(n)=-\frac{d}{b} s_2 x_2(n)^{-b}$.
\end{itemize}
\end{proof}


\begin{thebibliography}{99}


\bibitem{bahloulJ}
R. Bahloul,
\emph{Algorithm for computing Bernstein-Sato ideals associated with
a polynomial mapping},
J. Symbolic Comput.  32 (2001)  no. 6, 643--662.


\bibitem{bahloulC}
R. Bahloul,
\emph{D\'emonstration constructive de l'existence de polyn\^omes de
Bernstein-Sato pour plusieurs fonctions analytiques},
Compos. Math. 141 (2005), no. 1, 175--191.


\bibitem{bernstein}
I. N. Bernstein,
\emph{The analytic continuation of generalized functions with respect to a
parameter},
Funct. Anal. Appl. 6 (1972), 273--285.


\bibitem{bjork}
J.E. Bj\"ork,
\emph{Dimensions of modules over algebras of differential operators},
Fonctions analytiques de plusieurs variables et analyse complexe
(Colloq. Internat. CNRS, No. 208, Paris, 1972),  pp. 6--11.
``Agora Mathematica'', No. 1, Gauthier-Villars, Paris, 1974.


\bibitem{buchberger}
B. Buchberger,
\emph{Ein algorithmisches Kriterium f\"ur die L\"osbarkeit eines
algebraischen Gleichungssystems},
Aequationes Math. 4 (1970), 374--383.


\bibitem{briancon}
J. Brian\c{c}on,
\emph{Weierstrass pr\'epar\'e \`a la Hironaka},
Asterisque, Nos. 7 et 8 (1973), Soc. Math. France, 67--73.


\bibitem{bgmm}
J. Brian\c{c}on, M. Granger, Ph. Maisonobe, M. Miniconi,
\emph{Algorithme de calcul du polyn\^ome de Bernstein: cas non
d\'eg\'en\'er\'e},
Ann. Inst. Fourier (Grenoble) 39 (1989), no. 3, 553--610.


\bibitem{bmm1}
J. Brian\c{c}on, Ph. Maisonobe, M. Merle,
\emph{\'Eventails associ\'ees \`a des fonctions analytiques},
Tr. Mat. Inst. Steklova 238 (2002), Monodromiya v Zadachakh Algebr. Geom.
i Differ. Uravn., 70--80.


\bibitem{bmm2}
J. Brian\c{c}on, Ph. Maisonobe, M. Merle,
\emph{\'Equations fonctionnelles associ\'ees \`a des fonctions analytiques},
Tr. Mat. Inst. Steklova  238 (2002), Monodromiya v Zadachakh Algebr. Geom.
i Differ. Uravn., 86--96.


\bibitem{bmay}
J. Brian\c{c}on, H. Maynadier,
\emph{Equations fonctionnelles g\'en\'eralis\'ees: transversalit\'e et
principalit\'e de l'id\'eal de Bernstein-Sato},
J. Math. Kyoto Univ. 39 (1999), no. 2, 215--232.


\bibitem{cg}
F.J. Castro-Jim\'enez, M. Granger,
\emph{Explicit calculations in rings of differential operators},
S\'eminaires et Congr\`es 8 (2004), Soc. Math. France, 89--128.


\bibitem{kashiwara}
M. Kashiwara,
\emph{B-functions and holonomic systems. Rationality of roots of B-functions},
Invent. Math. 38 (1976/77), no. 1, 33--53.


\bibitem{kk}
M. Kashiwara, T. Kawai,
\emph{On holonomic systems for $\Pi \sp{N}\sb{l=1}(f\sb{l}+
(\surd -1)0)\sp{\lambda \sb{l}}$}, 
Publ. Res. Inst. Math. Sci. 15 (1979), no. 2, 551--575.


\bibitem{malgrange}
B. Malgrange,
\emph{Le polyn\^ome de Bernstein d'une singularit\'e isol\'ee},
Lecture Notes in Math. 459 (1975), Springer, Berlin , 98--119.


\bibitem{maynadierTez}
H. Maynadier,
\emph{Equations fonctionnelles pour une intersection compl\`ete
quasi-homog\`ene \`a singularit\'e isol\'ee et un germe semi-quasi-homog\`ene},
Th\`ese, Nice Sophia-Antipolis, 1996.


\bibitem{maynadierBull}
H. Maynadier,
\emph{Polyn\^omes de Bernstein-Sato associ\'es \`a une intersection
compl\`ete quasi-homog\`ene \`a singularit\'e isol\'ee},
Bull. Soc. Math. France 125 (1997), no. 4, 547--571.


\bibitem{oaku}
T. Oaku,
\emph{An algorithm of computing b-functions},
Duke Math. Journal 87 (1997), Vol 1, 115--132.


\bibitem{sabbah1}
C. Sabbah,
\emph{Proximit\'e \'evanescente I. La structure polaire d'un
$\mathcal{D}$-Module},
Compositio Math. 62 (1987), 283--328.


\bibitem{sabbah2}
C. Sabbah,
\emph{Proximit\'e \'evanescente II. Equations fonctionnelles pour plusieurs
fonctions analytiques},
Compositio Math. 64 (1987), 213--241.


\bibitem{torrelli}
T. Torrelli,
\emph{Polyn\^omes de Bernstein associ\'es \`a une fonction sur une
intersection compl\`ete \`a singularit\'e isol\'ee},
Ann. Inst. Fourier (Grenoble) 52 (2002), no. 1, 221--244.


\bibitem{ucha-castro}
J. M. Ucha, F. J. Castro-Jim\'enez,
\emph{On the computation of Bernstein-Sato ideals},
J. Symbolic Comput.  37 (2004), no. 5, 629--639.


\bibitem{walther}
U. Walther,
\emph{Bernstein-Sato polynomial versus cohomology of the Milnor fiber for
generic hyperplane arrangements},
Compositio Math. 141 (2005), no. 1, 121--145.



\end{thebibliography}
\end{document}